\documentclass{amsart}

\usepackage{amsmath}
\usepackage{amsthm}
\usepackage{amsopn}
\usepackage{amssymb}
\usepackage[cmtip, all]{xy}
\usepackage{stmaryrd}
 
\theoremstyle{plain}
\newtheorem{Lem}{Lemma}[section]

\newtheorem{Thm}[Lem]{Theorem}

{\theoremstyle{definition} 

}

\DeclareMathAlphabet{\mathpzc}{OT1}{pzc}{m}{it}
 
\DeclareMathOperator{\Map}{Map^\mathit{c}}

\DeclareMathOperator*{\holim}{holim}

\DeclareMathOperator*{\colim}{colim}

\newcommand{\zig}{\addtocounter{Lem}{1}\tag{\theLem}} 
 
\def\:{\colon}

\begin{document}
\title{Every $K(n)$-local spectrum is the homotopy fixed 
points of its Morava module}
\author{Daniel G. Davis$\sp 1$}
\author{Takeshi Torii}
\begin{abstract}
Let $n \geq 1$ and let $p$ be any prime. Also, let $E_n$ be the 
Lubin-Tate spectrum, $G_n$ the extended Morava stabilizer group, 
and $K(n)$ the $n$th Morava $K$-theory spectrum. Then work of 
Devinatz and Hopkins and some results due to Behrens and the first 
author of this note, show that if $X$ is a finite spectrum, then 
the localization $L_{K(n)}(X)$ is equivalent to the homotopy fixed point spectrum 
$(L_{K(n)}(E_n \wedge X))^{hG_n}$, which is formed with respect 
to the continuous action of $G_n$ on
$L_{K(n)}(E_n \wedge X)$. 
In this note, we show that 
this equivalence holds for {\em any} ($S$-cofibrant) spectrum $X$. Also, 
we show that for all such $X$, 
the strongly convergent Adams-type spectral sequence 
abutting to $\pi_\ast(L_{K(n)}(X))$ is isomorphic to the descent 
spectral sequence that abuts to 
$\pi_\ast((L_{K(n)}(E_n \wedge X))^{hG_n}).$
\end{abstract}

\footnotetext[1]{The first author was partially supported by a grant 
(\# $\negthinspace \negthinspace \negthinspace$ LEQSF(2008-11)-RD-A-27) from the Louisiana Board of Regents.} 

\maketitle
\par

\section {Introduction}
\par
In this note, we extend a result about the $K(n)$-localization 
of finite spectra, which is due to a combination of work by Devinatz and 
Hopkins (in \cite{DH}) and Behrens and the first 
author of this note (in \cite{joint}, \cite{cts}) (with most of the hard work being done 
by Devinatz and Hopkins), to all $K(n)$-local spectra.
\par
In more 
detail, let $n \geq 1$ and let $p$ be a prime. Above 
and elsewhere, $K(n)$ denotes the $n$th Morava $K$-theory spectrum, 
\[G_n = S_n \rtimes \mathrm{Gal}(\mathbb{F}_{p^n}/\mathbb{F}_p)\] 
is the $n$th extended Morava stabilizer group, and $E_n$ 
is the $n$th Lubin-Tate spectrum, with
\[\pi_\ast(E_n)=W(\mathbb{F}_{p^n})\llbracket u_1, ..., 
u_{n-1}\rrbracket[u^{\pm 1}],\]
where $W(\mathbb{F}_{p^n})$ is the ring of Witt vectors 
with coefficients in 
the field $\mathbb{F}_{p^n}$, the degree of $u$ is $-2$, 
and the complete power series ring 
is in degree zero. 
\par
Given a spectrum $X$, we define \[L_{K(n)}(E_n \wedge X),\] the 
Bousfield localization of $E_n \wedge X$ with respect to 
$K(n)$, to be the {\em point-set level Morava module} of $X$. 
When $\pi_\ast(L_{K(n)}(E_n \wedge X))$ satisfies certain hypotheses, 
it is common for these stable homotopy groups to be referred 
to as the {\em Morava module} of $X$. However, in this note, since 
we never use the term ``Morava module" in this algebraic sense, 
we will henceforth always refer to the point-set level Morava module 
of $X$ as just its Morava module.
\par
By \cite{Pgg/Hop0}, $G_n$ acts on $E_n$ through maps of commutative 
$S$-algebras and, by regarding $X$ as having trivial $G_n$-action and 
then giving its Morava module the diagonal $G_n$-action, 
$L_{K(n)}(E_n \wedge X)$ is a $G_n$-spectrum.
\par
Since the time of \cite{MRW} and the circulation of the results of 
\cite{Morava}, it has been believed 
by many experts in chromatic homotopy theory that it ought to 
be possible to realize the $K(n)$-localization of 
any spectrum $X$ as the $G_n$-homotopy 
fixed points of some $E_n$-module spectrum that is built out of 
$E_n$ and $X$. However, to date, it has not been clear how to 
do this. 
\par
Tremendous progress towards such a result was made 
by \cite{DH}, which showed that 
\begin{equation}\label{tremendous}\zig
E_n^{dhG_n} \simeq L_{K(n)}(S^0),\end{equation} where $E_n^{dhG_n}$ 
is a commutative $S$-algebra that behaves like a $G_n$-homotopy 
fixed point spectrum (e.g., the associated $K(n)$-local $E_n$-Adams spectral 
sequence looks like a descent spectral sequence, with $E_2$-term 
equal to continuous cohomology). 
\par
Notice that equivalence (\ref{tremendous}) implies that whenever 
$X$ is a finite spectrum,
\[E_n^{dhG_n} \wedge X \simeq L_{K(n)}(X).\] Additionally, by 
\cite[Theorem 1.3]{cts}, the Morava module 
$L_{K(n)}(E_n \wedge X) \simeq 
E_n \wedge X$ (since $X$ is finite) is a continuous $G_n$-spectrum, 
so that its $G_n$-homotopy fixed point spectrum 
$(L_{K(n)}(E_n \wedge X))^{hG_n}$ 
can be formed. Also, by 
\cite[Theorem 8.2.1]{joint}, there is an equivalence
\[E_n^{dhG_n} \simeq E_n^{hG_n}.\] Taken together, the preceding 
conclusions imply that $E_n^{hG_n} \wedge X \simeq L_{K(n)}(X)$, 
and hence, since $X$ is a finite spectrum,
\begin{equation}\label{finite}\zig
(L_{K(n)}(E_n \wedge X))^{hG_n} \simeq L_{K(n)}(X),\end{equation} by 
\cite[Theorem 9.9]{cts}. Therefore, the $K(n)$-localization of 
any finite spectrum can be realized as the $G_n$-homotopy fixed 
points of its Morava module.
\par
Now we are ready to explain how our last conclusion is generalized 
in this note. First, we remark that, from this point onward, we 
always work in the stable model category of symmetric spectra 
of simplicial sets or in its homotopy category. Thus, a ``spectrum" is a 
symmetric spectrum of simplicial sets and, when we 
work with, for example, commutative algebras, these objects 
are always to be understood as referring to commutative algebras in 
the setting of symmetric spectra.
\par
Recall that \cite[Theorem 1.3]{cts} shows that for any spectrum $X$, the Morava 
module $L_{K(n)}(E_n \wedge X)$ is a continuous $G_n$-spectrum, 
where, as before, 
$X$ is regarded as having the trivial $G_n$-action and the Morava 
module has the diagonal $G_n$-action. Then, in this note, 
we generalize equivalence (\ref{finite}) 
in the following way.
\begin{Thm}\label{main}
If $X$ is any $S$-cofibrant spectrum, then 
\[(L_{K(n)}(E_n \wedge X))^{hG_n} \simeq L_{K(n)}(X).\]
\end{Thm}
\par
We quickly make a technical (but useful) comment about 
Theorem \ref{main}. By using cofibrant replacement in the $S$ model structure on 
the category of symmetric spectra, given any spectrum $Z$, there 
is a weak equivalence $Z_c \rightarrow Z$ in the usual stable 
model category of symmetric spectra, with $Z_c$ $S$-cofibrant (see 
\cite[Section 5.3]{HSS}). Thus, there is no loss of generality in 
Theorem \ref{main} in requiring that $X$ be $S$-cofibrant, so that 
the theorem can be thought of as being valid for an arbitrary spectrum 
$X$. 
\par
Theorem \ref{main} shows that the $K(n)$-localization of any 
($S$-cofibrant) spectrum is the 
$G_n$-homotopy fixed points of its Morava module, answering the 
relatively old question of how to show that every $K(n)$-local spectrum can be 
obtained from a homotopy fixed point construction involving $E_n$ 
and $G_n$. The proof of Theorem \ref{main} is given in 
Section \ref{proof}.
\par
We give Theorem \ref{main} the desired (but generally unwieldy) 
``computational legs" with 
the following result.
\begin{Thm}\label{s.s.}
If $X$ is any $S$-cofibrant spectrum, then the strongly convergent $K(n)$-local 
$E_n$-Adams spectral sequence abutting to $\pi_\ast(L_{K(n)}(X))$ is 
isomorphic to the descent spectral sequence that abuts to 
$\pi_\ast((L_{K(n)}(E_n \wedge X))^{hG_n})$, from the $E_2$-terms 
onward. 
\end{Thm}
\par
If $X$ is a finite spectrum, then Theorem \ref{s.s.} reduces to 
\cite[Theorem 8.2.5]{joint}. We refer 
the reader to \cite[Appendix A]{DH} for an exposition of the aforementioned 
Adams-type spectral sequence. The descent spectral sequence that 
Theorem \ref{s.s.} refers to is defined in \cite[Section 4.6]{joint}. 
The proof of Theorem \ref{s.s.} is given in Section \ref{last}. 
\vspace{.1in}
\par 
\noindent
\textbf{Acknowledgements.} The proof of Theorem \ref{main} 
that appears in this 
note is a simplified version of an argument that relied more heavily 
on some results from \cite{joint}. Thus, the first author thanks Mark Behrens for 
various things that he learned from him during their 
collaboration on \cite{joint}. 
\section{The proof of Theorem \ref{main}}\label{proof}
\par
We begin this section by establishing some notation. 
We let \[c \: \widehat{S} \overset{\simeq}{\longrightarrow} 
L_{K(n)}(S^0)\] be a cofibrant replacement of 
$L_{K(n)}(S^0)$ in the model category of commutative symmetric ring 
spectra (see \cite[the discussion just before Theorem 19.6]{mmss}); 
the map $c$ is a weak equivalence in the stable model category of 
symmetric spectra. (We need the cofibrant commutative symmetric ring 
spectrum $\widehat{S}$ because later we will regard it as the ground 
ring of a profinite Galois extension.) Also, 
for the remainder of this section, we write $K$ 
in place of $K(n)$, so that our notation 
does not become too cumbersome. We will sometimes 
use the terminology of 
\cite[Section 1]{rocket}, adapted to the $K$-local category, 
as in \cite[Appendix A]{DH}.
\par
As in 
\cite[Section 5.2]{joint}, let $\mathpzc{Alg}$ be the model category 
of discrete commutative $G_n$-$\widehat{S}$-algebras: 
the objects of $\mathpzc{Alg}$ are 
discrete $G_n$-spectra that are also commutative $\widehat{S}$-algebras, 
and the morphisms are $G_n$-equivariant maps of commutative 
$\widehat{S}$-algebras. Let $(-)_\mathpzc{F} \: \mathpzc{Alg} \rightarrow 
\mathpzc{Alg}$ be a fibrant replacement functor 
for the model category $\mathpzc{Alg}$, 
and let $N \vartriangleleft_o G_n$ denote an 
open normal subgroup of $G_n$. 
Also, recall from \cite{DH} 
that each $E_n^{dhN}$ -- the commutative $\widehat{S}$-algebra that is 
written as $E_n^{hN}$ in \cite{DH} and behaves like the $N$-homotopy 
fixed point spectrum of $E_n$ -- 
is a $G_n/N$-spectrum that is $K$-local. 
Then, as in \cite{cts} and \cite{joint}, let 
\[F_n 
=  \colim_{N \vartriangleleft_o G_n} (E_n^{dhN})_\mathpzc{F};\] 
by construction, $F_n$ 
is a discrete $G_n$-spectrum and 
a commutative 
symmetric ring spectrum 
that is $E(n)$-local. Here, $E(n)$ is 
the usual Johnson-Wilson spectrum, with 
\[\pi_\ast(E(n)) = \mathbb{Z}_{(p)}[v_1, ..., v_{n-1}][v_n^{\pm1}].\]
\par
Let $X$ be any $S$-cofibrant spectrum. By \cite[Theorem 9.7]{cts}, 
\[(L_{K}(E_n \wedge X))^{hG_n} \simeq L_{K}((F_n \wedge X)^{hG_n}),\] 
where $(F_n \wedge X)^{hG_n}$ is the $G_n$-homotopy fixed points of the 
discrete $G_n$-spectrum $F_n \wedge X$. Thus, to prove Theorem \ref{main}, 
we only have to show that 
\[L_K((F_n \wedge X)^{hG_n}) \simeq L_{K}(X).\]
\par
Let $L_K\bigl((E_n)^{\wedge (\bullet+1)}\bigr)$ be the usual cosimplicial spectrum 
that is built from the unit map $S^0 \rightarrow 
E_n$ and the multiplication $E_n \wedge E_n \rightarrow E_n$. 
Here, for each $k \geq 0$,
\[L_K\bigl((E_n)^{\wedge{(k+1)}}\bigr) = 
L_K(\,\underbrace{E_n \wedge E_n \wedge 
\cdots \wedge E_n}_{k+1}\,)\] and there is the associated 
augmented resolution
\begin{equation}\label{resone}\zig
\ast \rightarrow \widehat{S} \rightarrow E_n \rightarrow 
L_K(E_n \wedge E_n) \rightarrow L_K(E_n \wedge E_n \wedge E_n) 
\rightarrow \cdots,\end{equation} which is the canonical 
$K$-local $E_n$-resolution of $\widehat{S}$. 
\par
Given $\widehat{S}$-modules $M$ and $N$, we let $M \wedge_{\widehat{S}}N$ 
denote their smash product in the category of $\widehat{S}$-modules. 
Then resolution (\ref{resone}) can be identified with the resolution 
\begin{equation}\label{restwo}\zig
\ast \rightarrow \widehat{S} \rightarrow E_n \rightarrow 
L_K(E_n \wedge_{\widehat{S}} E_n) \rightarrow 
L_K(E_n \wedge_{\widehat{S}} E_n \wedge_{\widehat{S}} E_n) 
\rightarrow \cdots,\end{equation} since, for each $k \geq 0$, 
\[L_K\bigl((E_n)^{\wedge{(k+1)}}\bigr) \simeq 
L_K\bigl((E_n)^{\wedge_{\widehat{S}}{(k+1)}}\bigr).
\] 
Since $E_n \simeq 
L_K(F_n)$ (this equivalence is due to \cite{DH}, but the reader might 
find the proof of it in \cite[Theorem 6.3]{cts} useful), 
resolution (\ref{restwo}) can be identified with the resolution 
\begin{equation}\label{resthree}\zig
\ast \rightarrow \widehat{S} \rightarrow L_K(F_n) \rightarrow 
L_K(F_n \wedge_{\widehat{S}} F_n) \rightarrow 
L_K(F_n \wedge_{\widehat{S}} F_n \wedge_{\widehat{S}} F_n) 
\rightarrow \cdots.\end{equation}
\par
For the next step in our proof, we make a few 
recollections. We define 
the cosimplicial spectrum $\mathrm{Map}^c(G_n^\bullet, F_n)$ as 
in \cite[Section 3.2]{joint}: if $K$ is a profinite group, then 
\[\mathrm{Map}^c(K, F_n) = \colim_{U \vartriangleleft_o K} 
\mathrm{Map}(K/U, F_n) \cong \colim_{U \vartriangleleft_o K} 
\textstyle{\prod}_{_{K/U}} F_n;\] $\mathrm{Map}^c(G_n, -)$ is a 
coaugmented comonad on the category of spectra; and, via the 
comonadic cobar construction, 
$\mathrm{Map}^c(G_n^\bullet, F_n)$ is the associated cosimplicial 
spectrum, which, in each codegree $k$, satisfies the isomorphism
\[\bigl(\mathrm{Map}^c(G_n^\bullet,F_n)\bigr)^k 
\cong \mathrm{Map}^c(G_n^k, F_n).\]
\par
We recall that by 
\cite[Theorem 5.4.4]{Rognes} and \cite[pg. 5034]{joint}, the work of 
\cite{DH} and \cite{LHS} implies that $F_n$ 
is a $K$-local profinite 
$G_n$-Galois extension of $\widehat{S}$ of finite vcd (in the sense 
of \cite{joint}). Thus, 
as in the proof of \cite[Proposition 6.3.1]{joint}, there is an 
equivalence
\[L_K((F_n)^{\wedge_{\widehat{S}}(\bullet +1)}) \simeq 
L_K(\Map(G_n^\bullet,F_n))\] 
of cosimplicial spectra. (Here, 
if a map $X^\bullet \rightarrow Y^\bullet$ of cosimplicial spectra 
is a weak equivalence in each 
codegree, then we regard it as a weak equivalence of cosimplicial 
spectra.) 
Therefore, resolution (\ref{resthree}) is 
equivalent 
to the resolution 
\begin{equation}\label{resfour}\zig
 \ast \rightarrow \widehat{S} \rightarrow L_K(F_n) \rightarrow 
L_K(\Map(G_n,F_n)) \rightarrow 
L_K(\Map(G_n^2,F_n)) 
\rightarrow \cdots,\end{equation} 
and hence, resolution (\ref{resfour}) can be regarded as 
the canonical $K$-local $E_n$-resolution of 
$\widehat{S}$.
\par
Below, to save space, we sometimes use 
the notation $Y \, \widehat{\wedge} \, Z$ 
to denote $L_K(Y \wedge Z)$, where $Y$ and $Z$ are 
arbitrary spectra, and, for the same reason, we occasionally write 
$(-)_K$ in place 
of $L_K(-)$. By smashing resolution (\ref{resfour}) 
with $X$ and then localizing with respect 
to $K$, it follows from 
\cite[Remark A.9]{DH} 
that 
\begin{equation*}
\ast \shortrightarrow X_K \shortrightarrow (F_n)_K 
\, \widehat{\wedge} \, X
\shortrightarrow \Map(G_n,F_n)_K \, \widehat{\wedge} \, X
\shortrightarrow 
\Map(G_n^2,F_n)_K \, \widehat{\wedge} \, X
\shortrightarrow \cdots\end{equation*} is a 
$K$-local $E_n$-resolution of 
$L_K(X)$. Thus, the equivalent resolution
\begin{equation*}\label{resfive}\zig
\ast \shortrightarrow X_K \shortrightarrow F_n 
\, \widehat{\wedge} \, X
\shortrightarrow \Map(G_n,F_n) \, \widehat{\wedge} \, X
\shortrightarrow 
\Map(G_n^2,F_n) \, \widehat{\wedge} \, X
\shortrightarrow \cdots\end{equation*} is a $K$-local $E_n$-resolution 
of $L_K(X)$, so that, since (\ref{resfive}) is the 
resolution associated to the canonical map 
\[ \Phi \: L_K(X) \rightarrow \holim_\Delta L_K(\Map(G_n^\bullet,F_n) \wedge 
X),\] $\Phi$ is a weak equivalence, 
by \cite[Corollary A.8]{DH}. 
\par
Since $F_n$ is $E(n)$-local and $L_{E(n)}(-)$ is a smashing localization, 
$F_n \wedge X$ 
is also $E(n)$-local. 
Then the proof of Theorem \ref{main} is finished by noting that
\begin{align*}
L_K((F_n \wedge X)^{hG_n}) & 
\simeq L_K(\holim_\Delta \Map(G_n^\bullet, F_n \wedge X)) \\
& \simeq \holim_\Delta L_K(\Map(G_n^\bullet, F_n \wedge X)) 
\\ & \simeq \holim_\Delta L_K(\Map(G_n^\bullet,F_n) \wedge X) 
\\
& \simeq L_K(X),
\end{align*}
where the first equivalence follows immediately from 
\cite[Theorem 3.2.1]{joint}; the second equivalence 
follows from the observation that 
$\Map(G_n^\bullet, F_n \wedge X)$ is a diagram of 
$E(n)$-local spectra (see, for example, \cite[Corollary 6.1.3; 
Lemma 6.1.4, $\negthinspace \negthinspace \negthinspace$ (3)]{joint}); and 
the third equivalence uses the fact that 
for each $k \geq 1$, there are natural equivalences 
\begin{align*}
\Map(G_n^k, F_n \wedge X) & 
\simeq \colim_{U \vartriangleleft_o G_n^k} \textstyle{\bigvee}_{G_n^k/U} 
(F_n \wedge X) \\ & \cong \colim_{U \vartriangleleft_o G_n^k} 
\bigl(\bigl(\textstyle{\bigvee}_{G_n^k/U}F_n\bigr) \wedge X\bigr)
\\ & \simeq \colim_{U \vartriangleleft_o G_n^k} 
\bigl(\bigl(\textstyle{\prod}_{_{G_n^k/U}}F_n\bigr) \wedge X\bigr)
\\ & \cong \Map(G_n^k, F_n) \wedge X,\end{align*} where the 
third equivalence (just above, between the two colimits) 
uses the fact that, since $X$ is $S$-cofibrant, 
the functor $(-) \wedge X$ preserves weak equivalences, 
by \cite[Corollary 5.3.10]{HSS}.
\section{The proof of Theorem \ref{s.s.}}\label{last}
\par
Let $X$ be any $S$-cofibrant spectrum. In 
Section \ref{proof}, we showed that resolution (\ref{resfive}) is a 
$K(n)$-local $E_n$-resolution of $L_{K(n)}(X)$. Thus, by 
\cite[discussion preceding Proposition A.5]{DH}, there is a map 
$\phi$ from the strongly convergent 
$K(n)$-local $E_n$-Adams spectral sequence 
for 
$\pi_\ast(L_{K(n)}(X))$, which we denote by $^\mathcal{A}E_r^{\ast, \ast}(X)$, 
to the homotopy spectral 
sequence for 
\[\pi_\ast\bigl(\holim_\Delta 
L_{K(n)}(\Map(G_n^\bullet,F_n) \wedge X)\bigr),\] which we denote by 
${^I}E_r^{\ast, \ast}(X)$. Furthermore, by \cite[Proposition A.5]{DH}, the 
map $\phi$ of 
spectral sequences is an isomorphism, from the $E_2$-terms onward.
\par
Notice that there are equivalences 
\begin{align*}
L_{K(n)}(\Map(G_n^\bullet,F_n) \wedge X) & \simeq 
L_{K(n)}(\Map(G_n^\bullet, F_n \wedge X)) \\
& \simeq \holim_i \Map(G_n^\bullet, F_n \wedge X \wedge M_i)\end{align*} 
of cosimplicial spectra, where the first equivalence was obtained at 
the end of Section~\ref{proof} and the second equivalence is an application of 
the fact that for any $E(n)$-local spectrum $Z$, $L_{K(n)}(Z) 
\simeq \holim_i (Z \wedge M_i)$, where
\[M_0 \leftarrow M_1 \leftarrow \cdots \leftarrow M_i \leftarrow \cdots \] is a suitable 
tower of generalized Moore spectra (see \cite[Section 2]{HoveyCech}). 
\par
The above equivalences of cosimplicial spectra 
imply that ${^I}E_r^{\ast, \ast}(X)$ is isomorphic to the 
homotopy spectral sequence for 
\[\pi_\ast\bigl(\holim_\Delta \holim_i \Map(G_n^\bullet, 
F_n \wedge X \wedge M_i)\bigr),\] which we denote by ${^{II}}E_r^{\ast, \ast}(X)$. 
By \cite[Section 4.6]{joint}, ${^{II}}E_r^{\ast, \ast}(X)$ is exactly the 
descent spectral sequence 
for $\pi_\ast((L_{K(n)}(E_n \wedge X))^{hG_n})$, so that the 
chain of isomorphisms 
\[^\mathcal{A}E_r^{\ast, \ast}(X) 
\cong {^{I}}E_r^{\ast, \ast}(X)  \cong {^{II}}E_r^{\ast, \ast}(X)\]
completes the proof of Theorem \ref{s.s.}.

\end{document}